\documentclass{article}

\usepackage{arxiv}

\usepackage[utf8]{inputenc} 
\usepackage[T1]{fontenc}    
\usepackage{hyperref}       
\usepackage{url}            
\usepackage{booktabs}       
\usepackage{amsfonts}       
\usepackage{nicefrac}       
\usepackage{microtype}      
\usepackage{lipsum}

\RequirePackage{amsthm,amsmath}

\usepackage{graphicx}
\usepackage{amssymb,bm}
\usepackage{verbatim}
\usepackage{multirow}
\usepackage{color,soul}
\usepackage[dvipsnames]{xcolor}

\newtheorem{remark}{Remark}

\title{Reply to ``On some problems in the article Efficient likelihood estimation in state space models''}

\author{
  Cheng-Der Fuh
  \\
  Fanhai International School of Finance\\
  Fudan University\\
  Shanghai 200433, China \\
  \texttt{cdffuh@gmail.com} \\
   \And
 Chu-Lan Kao \\
  Institute of Statistics\\
  National Chiao-Tung University\\
  Hsinchu 30010, Taiwan \\
  \texttt{chulankao@gmail.com} \\
}

\begin{document}
\maketitle



\def\theequation{1.\arabic{equation}}
\setcounter{equation}{0}

The first author, Cheng-Der Fuh, is grateful for the comments by Dr. Jensen (2010). This note replies his comments on Problem 2.3, which was left in Fuh (2010). In the following, we use the same notations and definitions in Fuh (2006) unless specified.\\

\noindent
{\bf Problem 2.3 Asymptotic properties of score function and observed information}

Page 2060, L12.
In the proof of Lemma 6, (7.9) defined a new iterated functions system, therefore Corollary 1 can not be used directly. 
The same situation happens for Theorems 5 and 7. The rigorous proofs of these results will be given in a separate paper.\\~

Let ${\bf X}= \{X_n, n \geq 0 \}$ be a Markov
chain on a general state space ${\cal X}$, with transition
probability kernel $P^{\theta}(x,\cdot)= P^{\theta}\{X_1 \in \cdot|X_0=x\}$ and
stationary probability $\pi(\cdot):=\pi_{\theta}(\cdot)$, where $\theta \in \Theta
\subseteq {\bf R}^q $ denotes the unknown parameter. Let $\xi_{1:n}$ be the observation such that the distribution of $\xi_n$ depends on $X_n$ and $\xi_{n-1}$, independent to others. Denote $L(\theta;\xi_{1:n})$ as the full likelihood of $\xi_{1:n}$, and $\ell(\theta)= \log L(\theta;\xi_{1:n})$ as the $\log$ likelihood. Fuh (2006) has shown that $\ell(\theta)$ can be written as an additive functional of a Markovian iterated random function system (MIRFS). In the following, we will show that the derivatives of $\ell(\theta)$ can also be written as an additive functional of a MIRFS.

To start with, let
\begin{eqnarray}\label{mn}
M_n:= {\bf P}_\theta(\xi_n) \circ \cdots \circ {\bf P}_\theta(\xi_1) \circ {\bf P}_\theta(\xi_0)
\end{eqnarray}
be defined in (5.6) of Fuh (2006). Then, as in Fuh (2006), (5.7), we have
\begin{eqnarray}\label{logLik}
\ell(\theta) := \log L(\theta;\xi_{1:n}) = \log \Vert M_n \pi \Vert= \sum_{i=1}^n \log \frac{\Vert M_n \pi \Vert}{\Vert M_{n-1} \pi \Vert}
= \sum_{i=1}^n g^0( M_n^0, M_{n-1}^0). 
\end{eqnarray}
Hence, $\ell(\theta)$ is an additive functional of  $((X_n, \xi_n), M_n)$.
Fuh (2006) shows that $\{((X_n, \xi_n), M_n), n \geq 0\}$ forms an ergodic Markov chain, induced by the MIRFS based on (\ref{mn}), 
on the state space $(\mathcal{X} \times {\bf R}^d) \times {\bf M}$, where ${\bf M}$ is the function space defined in Fuh (2006), page 2046.  Then Fuh (2006) uses this result to prove the law of large number for the log likelihood.
In this note, we will show that this result can be extended to higher derivatives of the log likelihood.

To this end, for any $1 \leq i \leq q$ and positive integer $k$, let $D_i$ be the partial derivative with respective to the $i$-th dimension of $\theta$ in some neighborhood $N_\delta(\theta_0) := \{\theta: |\theta - \theta_0| < \delta\}$ of the true value $\theta_0$, and $(D_i)^k$ be the corresponding $k$th partial derivative. Note that for any two given random functions ${\bf P}_\theta(\xi_{j+1})$ and ${\bf P}_\theta(\xi_{j})$, and any $h_\theta$, by  Assumptions C1, C2'-C5', and C6-C9 in  Fuh (2006)  and the dominated convergence  theorem, we have
	\begin{eqnarray*}
	&~& D_i \left\{ {\bf P}_\theta(\xi_{j}) h_\theta(x) \right\}
	= D_i \left\{ 
	\int_{y \in \mathcal{X}} p_\theta(y,x) f(\xi_j; \theta|x, \xi_{j-1})h_\theta (y)m(dy) \right\} \\
	&= &
	\int_{y \in \mathcal{X}} \bigg\{ f(\xi_j; \theta|x, \xi_{j-1}) h_\theta(y) D_i p_\theta(y,x) + p_\theta(y,x)  h_\theta(y) D_i f(\xi_j; \theta|x, \xi_{j-1})  \\
	&~&~~~~~~~~~ 
	+ p_\theta(y,x) f(\xi_j; \theta|x, \xi_{j-1}) D_i h_\theta(y) \bigg\} m(dy),
	\end{eqnarray*}
	and
	\begin{eqnarray*}
	&~& D_i \left\{ {\bf P}_\theta(\xi_{j+1}) \circ {\bf P}_\theta(\xi_{j}) h_\theta(x) \right\} \\
	&=& 	D_i \bigg\{ \int_{z \in \mathcal{X}} p_\theta(z,x) f(\xi_{j+1}; \theta|x, \xi_j) 
    \times \left( 
	\int_{y \in \mathcal{X}} p_\theta(y,z) f(\xi_j; \theta|z, \xi_{j-1})h_\theta(y)m(dy)
	\right) m(dz)  \bigg\} \\
	&=&
	\int_{z \in \mathcal{X}} D_i \left\{ p_\theta(z,x) f(\xi_{j+1}; \theta|x, \xi_j) \right\}
	\times \left( 
	\int_{y \in \mathcal{X}} p_\theta(y,z) f(\xi_j; \theta|z, \xi_{j-1})h_\theta(y)m(dy)
	\right) m(dz) \\
	&~&+
	\int_{z \in \mathcal{X}} p_\theta(z,x) f(\xi_{j+1}; \theta|x, \xi_j)
	\times \left( 
	\int_{y \in \mathcal{X}} D_i \left\{ p_\theta(y,z) f(\xi_j; \theta|z, \xi_{j-1})  h_\theta(y) \right\} m(dy)\right) m(dz) \\
	&=& \left\{ D_i {\bf P}_\theta(\xi_{j+1}) \right\} \circ {\bf P}_\theta(\xi_{j}) h_\theta(x) 
	+	{\bf P}_\theta(\xi_{j+1}) \circ \left\{ D_i ({\bf P}_\theta(\xi_{j}) h_\theta (x) )\right\}.  
	\end{eqnarray*}

For a given non-negative integer vector $\nu=(\nu^{(1)},\cdots,\nu^{(q)})$ write
$|\nu|=\nu^{(1)}+\cdots+\nu^{(q)}$, $\nu!=\nu^{(1)}! \cdots \nu^{(q)}!$,
and let $D^{\nu}=(D_1)^{\nu^{(1)}} \cdots (D_q)^{\nu^{(q)}}$ denote the
$\nu$th derivative with respect to $\theta$ in $N_\delta(\theta_0)$. For any $\nu$, define
\begin{eqnarray*}
    W_n^\nu = D^\nu M_n = 
     (D_1)^{\nu^{(1)}} \cdots (D_q)^{\nu^{(q)}} (M_n). 
\end{eqnarray*}
Then by  Assumptions C1, C2'-C5', and C6-C9 in  Fuh (2006)  and the dominated convergence
theorem, we have $ D^\nu \Vert (M_n \pi) \Vert = \Vert D^\nu (M_n \pi) \Vert$. 

Now let us consider all derivatives with order $r$ or less. Note that for a fixed integer $r \geq 1$,  there are exactly $K = {(r+q)!}/{r!q!}$ different $\nu$ satisfying $|\nu| \leq r$. Label all such $\nu$ by $\nu_1, \nu_2, \cdots, \nu_K$, and let $W_n = (W_n^{\nu_1}, W_n^{\nu_2}, \cdots, W_n^{\nu_K})^t$, with $t$ denotes the transpose. Then $W_n \in {\bf M}^K := \{ v=(m_1, \cdots, m_K)^t: m_k \in {\bf M}, 1 \leq k \leq K \}$. Moreover, for given $\nu_i$ and $\nu_j$, let $\nu_i + \nu_j$ denote addition of each component in the vector. 

To investigate the dynamic of $W_n$, note that for any $\nu_i$, we have
\begin{align}\label{iteration}
 W_n^{\nu_i} 
=&  D^{\nu_i} \big({\bf P}_\theta(\xi_n) \circ \cdots \circ {\bf P}_\theta(\xi_1) \circ {\bf P}_\theta(\xi_0) \big) \\
\notag
 =&  \sum_{\substack{1 \leq j \leq k \leq K \\ \nu_i = \nu_j + \nu_k}} \bigg\{ \frac{(\nu_i)!}{(\nu_j)! (\nu_k)!} D^{\nu_k} {\bf P}_\theta(\xi_n) \circ D^{\nu_j} \bigg( {\bf P}_\theta(\xi_{n-1}) \circ \cdots  \circ {\bf P}_\theta(\xi_1) \circ {\bf P}_\theta(\xi_0) \bigg) \bigg\} \\
\notag
 = & \sum_{\substack{1 \leq j \leq k \leq K \\ \nu_i = \nu_j + \nu_k}} \frac{(\nu_i)!}{(\nu_j)! (\nu_k)!} \left\{ D^{\nu_k} {\bf P}_\theta(\xi_n) \circ W_{n-1}^{\nu_j}  \right\}. 
\end{align}
Hence, we can define a $K$-by-$K$ {\it matrix form} $A_n = \left\{ a_n^{ij}: 1 \leq i, j \leq K \right\}$, with each $a_n^{ij} \in {\bf M}$ defined as
\begin{equation}\label{aij}
a_n^{ij} = \left\{
\begin{array}{cc}
\frac{(\nu_i)!}{(\nu_j)!(\nu_k)!} D^{\nu_k} {\bf P}_\theta(\xi_n) & \mbox{if }  \exists 1 \leq k \leq K \mbox{ such~that } \nu_i = \nu_j + \nu_k, \\
0 & \mbox{otherwise}.
\end{array}
\right.
\end{equation}

In addition, for each $K$-by-$K$ ${\bf M}$-valued matrix form $B = \{b_{ij}: 1 \leq i, j \leq K\}$, and each $K$-dimensional ${\bf M}$-valued vector $V = (V_1, V_2, \cdots, V_K) \in {\bf M}^K$,  we define
\begin{equation}\label{operator}
    B \circ V := \left(
\begin{array}{c}
\sum_{j=1}^K b_{1j} \circ V_j \\
\sum_{j=1}^K b_{2j} \circ V_j \\
\vdots\\
\sum_{j=1}^K b_{Kj} \circ V_j
\end{array}
\right).
\end{equation}
Then by (\ref{iteration}), we have $W_n =A_n \circ W_{n-1}$, and thus
\begin{eqnarray}\label{wn}
W_n = A_n \circ A_{n-1} \cdots A_1 \circ W_0,
\end{eqnarray}
where $W_0 = \{ W_0^\nu: |\nu| \leq r \}$ with $W_0^\nu =D^\nu {\bf P}_\theta(\xi_0) .$

\begin{remark}
To illustrate {\rm (\ref{wn})}, let $q=1$, i.e., $\theta$ is one dimensional. In this case, $\nu \in {\bf R}^1$ and we can simply label all $|\nu|\leq r$ by natural order so that $W_n = (W_n^0, W_n^1, \cdots, W_n^r)^t$, the vector of the first $r$-th derivatives. Then for any $0 \leq k \leq r$, we have
\begin{eqnarray*}
 W_n^k 
 & =& D^k ({\bf P}_\theta(\xi_n) \circ \cdots \circ {\bf P}_\theta(\xi_1) \circ {\bf P}_\theta(\xi_0)) \\
& =& \sum_{0 \leq k_1 \leq k} \bigg\{ \frac{k!}{(k_1)! (k-k_1)!} D^{k_1} 
{\bf P}_\theta(\xi_n) \circ D^{k-k_1} \bigg({\bf P}_\theta(\xi_{n-1}) \circ \cdots \circ {\bf P}_\theta(\xi_1) \circ {\bf P}_\theta(\xi_0) \bigg)\bigg\} \\
& = &  \sum_{0 \leq k_1 \leq k} C_{k_1}^k \left\{ D^{k_1} {\bf P}_\theta(\xi_n) \circ W_{n-1}^{k-k_1}  \right\},
\end{eqnarray*}
where $C_a^b = \frac{b!}{a!(b-a)!}$. Therefore  $W_n = A_n \circ W_{n-1}$ with
\begin{equation}
~~~~~~A_n = 
    \left(
    \begin{array}{cccc}
     {\bf P}_\theta(\xi_n) & 0 & \cdots & 0 \\ 
     C_1^1 D^1 {\bf P}_\theta(\xi_n) & {\bf P}_\theta(\xi_n) & \cdots & 0 \\
     \vdots & \vdots & \ddots & \vdots \\
     C_r^r D^r {\bf P}_\theta(\xi_n) & C_{r-1}^r D^{r-1} {\bf P}_\theta(\xi_n) & \cdots & {\bf P}_\theta(\xi_n)
    \end{array}
    \right), \label{An}
\end{equation}
where $0$ denotes the zero function in ${\bf M}$.
\end{remark}

\begin{remark}
Note that $A_n$  in {\rm (\ref{aij})} and $W_n$ in {\rm (\ref{wn})} are ${\bf M}$-valued, other than the traditional ${\bf R}$-valued vector and matrix, respectively. To illustrate this phenomenon, we consider a finite $D$-state Markov chain and an one-dimensional parameter $\theta$ case, then $A_n$ in {\rm (\ref{An})} is a $K$-by-$K$ matrix form with each element being a $D$-by-$D$ matrix (with $0$ being a $D$-by-$D$ zero matrix.) In the same matter, although the operator defined in {\rm (\ref{operator})} looks like a traditional matrix multiplication, it is different by having the multiplication within each component replaced by $\circ$.
Nevertheless, the essential idea is to introduce a matrix form for $W_n$, which can be used to show that it  forms an ergodic Markov chain via {\rm (\ref{wn})}.
\end{remark}

It is worth mentioning that the feature of getting a neat form in (\ref{wn}) is based on a matrix representation (\ref{aij}) for all partial derivatives up to the $r$th-order. Next by (\ref{wn}) and under the Assumptions C1, C2'-C5' and C6-C9 in  Fuh (2006), with the following additional condition:
\begin{itemize}
\item[C10] For any $\theta \in N_\delta(\theta_0)$ and $\nu$ with $|\nu| \leq r$, $\sup_{x \in \mathcal{X}} \vert D^\nu f(s_1; \theta|x, s_0) \vert < \infty$.
\end{itemize}
Then, through a process similar to the proof of Lemma 3 in Fuh (2006), it is straightforward to show  that for $\theta \in N_\delta(\theta_0)$, the MIRFS 
$\{((X_n, \xi_n), W_n), n \geq 0\}$ in (\ref{wn}) satisfies Assumption K in Fuh (2006).
Furthermore, Lemma 4 in Fuh (2006) holds for the induced Markov chain $\{((X_n, \xi_n), W_n), n \geq 0\}$ on the state space $({\cal X} \times {\bf R}^d) \times {\bf M}^K$. 

Last, we present a specific form of the first and second partial derivatives of the log likelihood function as follows. For $|\nu|= 1$, we have
\begin{align}
\label{1stderiv}
	& D^\nu (\log L( \theta; \xi_{1:n} ) ) = D^\nu ( \log \Vert M_n \pi \Vert)= 
	\frac{ \Vert D^\nu (M_n \pi )\Vert}{\Vert M_n \pi \Vert}\\
\notag
	= &
	\frac{ \Vert (D^\nu M_n) \pi + M_n(D^\nu \pi) \Vert}{\Vert M_n \pi\Vert}
	 \\
\notag
	 = &
	\sum_{i=1}^n \frac{ \Vert (W_i^\nu) \pi + W_i^0(D^\nu \pi) \Vert}{\Vert W_i^0 \pi\Vert} - 
\frac{ \Vert (W_{i-1}^\nu) \pi + W_{i-1}^0(D^\nu \pi) \Vert}{\Vert W_{i-1}^0 \pi\Vert}
	\\
\notag
	:= & \sum_{i=1}^n g^\nu ( W_i, W_{i-1}).
\end{align}
That is, the first derivative of the log likelihood function can be rewritten as an additive functional of the Markov chain $\{((X_n, \xi_n), W_n), n \geq 0\}$.

To represent the second partial derivative of the log likelihood, for $|\nu|=2$, let us write $\nu=\nu_1+\nu_2$ such that $|\nu_1| = |\nu_2|=1$. Then, we have
\begin{align*}
	& D^{\nu} \log L(\theta;\xi_{1:n}) = D^{\nu_1} \left( D^{\nu_2} \log L( \theta ; \xi_{1:n}) \right)
	= D^{\nu_1} \frac{\Vert W_n^{\nu_2} \pi + W_n^0 (D^{\nu_2} \pi)\Vert}{\Vert W_n^0 \pi \Vert} \\
	= & \frac{\Vert  W_n^{\nu} \pi + W_n^{\nu_2} (D^{\nu_1} \pi) + W_n^{\nu_1} (D^{\nu_2}\pi) + W_n^0 (D^{\nu}\pi)  \Vert}{\Vert W_n^0 \pi \Vert}  - \frac{\Vert W_n^{\nu_2} \pi + W_n^0 (D^{\nu_2} \pi)\Vert \times \Vert W_n^{\nu_1} \pi + W_n^0 (D^{\nu_1} \pi)\Vert}{\Vert W_n^0 \pi \Vert^2} \\
	:= & \sum g^\nu(W_i, W_{i-1})
\end{align*}
where 
\begin{align*}
& g^\nu(W_i, W_{i-1}) \\
= & 
\left\{ \frac{\Vert  W_i^{\nu} \pi + W_i^{\nu_2} (D^{\nu_1} \pi) + W_i^{\nu_1} (D^{\nu_2}\pi) + W_i^0 (D^{\nu}\pi)  \Vert}{\Vert W_i^0 \pi \Vert}  -
\frac{\Vert  W_{i-1}^{\nu} \pi + W_{i-1}^{\nu_2} (D^{\nu_1} \pi) + W_{i-1}^{\nu_1} (D^{\nu_2}\pi) + W_{i-1}^0 (D^{\nu}\pi)  \Vert}{\Vert W_{i-1}^0 \pi \Vert}
\right\}
\\
& -
\left\{ \frac{\Vert W_i^{\nu_2} \pi + W_i^0 (D^{\nu_2} \pi)\Vert \times \Vert W_i^{\nu_1} \pi + W_i^0 (D^{\nu_1} \pi)\Vert}{\Vert W_i^0 \pi \Vert^2}  - \frac{\Vert W_{i-1}^{\nu_2} \pi + W_{i-1}^0 (D^{\nu_2} \pi)\Vert \times \Vert W_{i-1}^{\nu_1} \pi + W_{i-1}^0 (D^{\nu_1} \pi) \Vert}{\Vert W_{i-1}^0 \pi \Vert^2} \right\}.
\end{align*}
That is, the second derivative of the log likelihood function can be rewritten as an additive functional of the Markov chain $\{((X_n, \xi_n), W_n), n \geq 0\}$.

In fact, for any $\nu$ with $|\nu| \leq r$, there exists a function $g^\nu$ such that
\begin{eqnarray*}
	D^\nu \log L(\theta;\xi_{1:n}) = \sum_{i=1}^n g^\nu (W_i, W_{i-1}).
\end{eqnarray*}
In other words, the $\nu$-th partial derivatives of the log likelihood can be rewritten as an additive functional of the Markov chain  $\{((X_n, \xi_n), W_n), n \geq 0\}$. This can be proved as follows. As stated in (\ref{1stderiv}), such $g^\nu$ exists for all $|\nu|=1$. Now for all $|\nu|$ with $|\nu| \in [2,r]$, take $\nu_1$ and $\nu_2$ such that $|\nu_1|=1$ and $\nu_1+\nu_2=\nu$. Then, we have
\begin{eqnarray*}
	D^{\nu_1+\nu_2} \log L(\theta;\xi_{1:n}) = D^{\nu_2} \sum_{i=1}^n g^{\nu_1} (W_i, W_{i-1}) = \sum_{i=1}^n (D^{\nu_2} g^{\nu_1} (W_i, W_{i-1})).
\end{eqnarray*}
Moreover, as shown in (\ref{1stderiv}), $g^{\nu_1} (W_i, W_{i-1})$ only involves the derivatives up to the order of $1$, so $(D^{\nu_2} g^{\nu_1} (W_i, W_{i-1}))$ only involves the derivatives up to the order of $|\nu|$, so it is still a function of $W_i$ and $W_{i-1}$ as they consist of all derivatives up to the order of $r$. Thus, such $g^{\nu}$ exists for all $|\nu| \leq r$.\\

\baselineskip = 16pt
\centerline{\bf REFERENCES}
\begin{description}


\item{} Fuh, C. D. (2006). Efficient likelihood estimation in state space models. {\it Ann. Statist.} 
{\bf 34}, 2026-2068.
\item{} Fuh, C. D. (2010). Reply to ``on some problems in the article efficient likelihood estimation in state space models". {\it Ann. Statist.} {\bf 38}, 1282-1285.
\item{} Jensen, J. L. (2010). On some problems in the article efficient likelihood estimation in state space models. {\it Ann. Statist.} {\bf 38}, 1279-1281.
\end{description}


\end{document}